\documentclass[12pt,oneside]{amsart}
\usepackage{geometry}                % See geometry.pdf to learn the layout options. There are lots.
\usepackage{graphicx}
\usepackage{tikz}
\usetikzlibrary{shapes,positioning,intersections,quotes}
\usepackage{amssymb}
\usepackage{epstopdf}

\usepackage{amsmath,amsthm,amscd,amssymb, enumerate}
\usepackage{latexsym}
\usepackage[colorlinks,citecolor=red,pagebackref,hypertexnames=false]{hyperref}
\numberwithin{equation}{section}
\theoremstyle{plain}
\newtheorem{theorem}{Theorem}[section]
\newtheorem{lemma}[theorem]{Lemma}
\newtheorem{corollary}[theorem]{Corollary}

\theoremstyle{definition}

\theoremstyle{remark}

\newtheorem{case[theorem]}{Case}

\def\bc{\begin{corollary}}
\def\ec{\end{corollary}}
\def\be{\begin{equation}}
\def\ee{\end{equation}}
\def\bast{\begin{eqnarray*} }
\def\east{\end{eqnarray*} }
\def\bea{\begin{eqnarray}}
\def\eea{\end{eqnarray}}

\def\bmat{\begin{matrix}}
\def\emat{\end{matrix}}

\def\R{\mathbb R}

\def\R{\Bbb R}

\def\({\left(}
\def\){\right)}
\def\[{\left[}
\def\]{\right]}
\def\<{\left\langle}
\def\>{\right\rangle}

\title{generalized Radon transforms on fractal measures}

\author{Shengze Duan}

\begin{document}
\maketitle

\begin{abstract}In the setting of a general Borel measure $\mu$ on $\R^d$ with the natural ball size condition $$\mu[B(x,r)]\leq Cr^s,$$ we establish the $L^p(\mu)$-$L^q(\mu)$-estimate for the generalized Radon transform $$(Af)(x):=\int_{\Phi(x,y)=0}(f\mu)(y)\psi(x,y)d\sigma_x(y),$$ where $\Phi$ is a smooth function away from the diagonal. 

Among other reasonable assumptions, an $L^2$-Sobolev bound on $A$ on $\R^d$ is imposed. This bound is satisfied in many natural situations. The main result is, in general, sharp up to endpoints. 
\end{abstract}  

\maketitle

%%%%%%%%%%%%%%%%%%%%

\section{Introduction} 
It was proved by Littman \cite{littman} and Strichartz \cite{strichartz} that the classical spherical averaging operator 
$$(Af)(x):=\int f(x-y)d\sigma(y),$$
where $\sigma$ denotes the surface measure on the unit sphere $\mathbb S^{d-1}$ ($d>1$ through this paper), is bounded from $L^p(\mathbb R^d)$ to $L^q(\mathbb R^d)$ under the Lebesgue Measure, if and only if $(\frac{1}{p},\frac{1}{q})$ belongs to the closed triangle region in the diagram
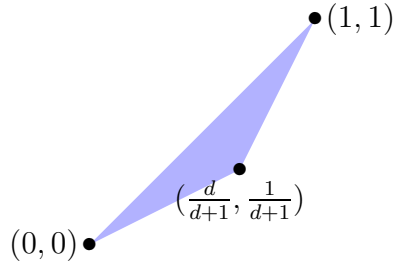
\begin{center}
\begin{tikzpicture}
\path[fill=blue!30] (0,0)--(3,3)--(2,1)--(0,0);
\draw (2,1) node[below]{$(\frac{d}{d+1},\frac{1}{d+1})$};
\draw (2,1) node{$\bullet$};
\draw (3,3) node[right]{$(1,1)$};
\draw (3,3) node{$\bullet$};
\draw (0,0) node[left]{$(0,0)$};
\draw (0,0) node{$\bullet$};
\end{tikzpicture}
\end{center}
As has been studied by Guillemin and Sternberg \cite{GS} and Phong and Stein \cite{PS}, the above operator is a special case of the generalized Radon transform
\be\label{radon}
(Af)(x):=\int_{\Phi(x,y)=0}f(y)\psi(x,y)d\sigma_x(y),
\ee
where $\sigma_x(y)$ is the surface measure on the submanifold defined by the zero set of $\Phi$ for every fixed $x$, and reasonable conditions on $\Phi$ can be imposed.

\smallskip
In this paper, we consider a compactly supported Borel measure $\mu$ satisfying
\be\label{ball}
\mu[B(x,r)]\leq Cr^s,
\ee
where $C$ and $s$ is uniform in all $x\in\R^d$ and all $r$ smaller than some fixed value and consider the operator \eqref{radon}
from $L^p(\mu)$ to $L^q(\mu)$. One can easily replace the codomain's measure by a different $\mu'$ with a different index $s'$, but for simplicity we set only one measure. The justification for $Af$, where $f\in L^1(\mu)$, is as follows. Fix a smooth non-negative compactly supported function $\phi$ with integral $1$, and use the approximation of the identity 
$$\phi_\epsilon(\cdot):=\frac{1}{\epsilon^d}\phi(\cdot)$$
to define $Af$ as the $L^q(\mu)$-limit of $A(\phi_\epsilon*f\mu)$. One may check the Cauchyness of the sequence via essentially the same proof as our argument. The function $\Phi(x,y)$ is assumed to be smooth on $(\R^d)^2$, and its gradients, $\nabla_x\Phi$ and $\nabla_y\Phi$, with respect to both $R^d$-variables $x,y$, are assumed to be nonzero on the support of some fixed smooth compactly supported $\psi$. Our main result is the following. 
\begin{theorem}\label{main}
Assume the condition \eqref{ball} and that both $A$ and its adjoint are bounded from $L^2(\R^d)$ to the $L^2$-based Sobolev space $L^2_\alpha(\R^d)$. Then, if $s>d-\alpha$, $A$ is bounded from $L^p(\mu)$ to $L^{p'}(\mu)$ for $\frac{1}{p}\in[\frac{1}{2},\frac{2s+2\alpha-2d+1}{2s+2\alpha-2d+2})$, and from $L^p(\mu)$ to $L^{2}(\mu)$ for $\frac{1}{p}\in[\frac{1}{2},\frac{3s+2\alpha-2d}{2s})$.

\smallskip
Moreover, if $\mu$ further satisfies the lower bound $\mu (B(x,r))\geq c r^s$ for $x$ on its support, then $A$ is bounded from $L^p(\mu)$ to $L^p(\mu)$ for $\frac{1}{p}\in[\frac{1}{2},\frac{s+2\alpha-d}{2\alpha})$. In general, exact endpoint values aside, the threshold for $s$ is sharp, and the Riesz diagram is complete with the $L^p$-$L^p$ and the $L^p$-$L^{p'}$-endpoint.
\end{theorem}\noindent

A few remarks are warranted here to clarify the nature of the result. If $\Phi$ has codomain $\R^n$, so the surface we integrate over in \ref{radon} is codimension-$n$, then the endpoints will reflect $n$ with the same proof under some appropriate nonvanishing condition on $\Phi$'s first derivatives. The Sobolev condition on the adjoint operator implies that our Riesz diagram will be symmetrical around the axis $p=q'$, and the $L^p$-$L^2$-endpoint is the only one that needs this condition, without which this endpoint is only on the symmetrical position around the axis. The compact support tells that one can always go upwards vertically in the diagram.
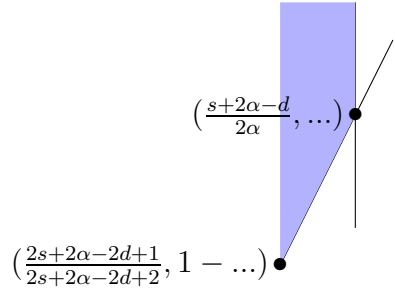
\begin{center}
\begin{tikzpicture}
\draw (0,0)--(1.5,3);
\draw (1,0.5)--(1,3.5);
\path[fill=blue!30] (0,0)--(1,2)--(1,3.5)--(0,3.5)--(0,0);
\draw (0,0) node[left]{$(\frac{2s+2\alpha-2d+1}{2s+2\alpha-2d+2},1-...)$};
\draw (0,0) node{$\bullet$};
\draw (1,2) node[left]{$(\frac{s+2\alpha-d}{2\alpha},...)$};
\draw (1,2) node{$\bullet$};
\end{tikzpicture}
\end{center}

\noindent
The $L^p$-$L^2$-endpoint, in general, can be on either side of or on the sloped sharpness line in the diagram. When $A$ is the spherical averaging operator, or $\Phi(x,y):=x\cdot y-1$ with $\psi$ supported away from the origin, $\alpha=\frac{d-1}{2}$. This is due to the nonzero rotational curvature. See, for example, Stein and Shakarchi \cite{stein} for reference. This $\alpha$ makes the three endpoints on the same line, and examples as such provide sharpness in dimension $2$ and $3$. But in dimension higher than $3$, to obtain the sharpness lines in the diagram, we need to set $\alpha\leq1$.

\smallskip
As for the threshold $s>d-\alpha$, similarly, in dimension higher than $3$, we need to set $\alpha\leq1$ to get sharpness. However, when $\alpha=\frac{d-1}{2}$ and $d>3$, we have a sharpness example of a parabolic hull convolution on a parabolic discrete measure in a uniform sense. That is, we allow the measure to vary but with uniform total mass and ball size condition \eqref{ball}, so in the corresponding positive result, operator norm should depend on the assumed uniform constants. The parabolic example is in Iosevich and Senger \cite{alex2}, and the reference to the original construction by Valtr, along with the connection to the Falconer distance problem, can be found there.

\smallskip 
For comparison purpose and as an $L^p$-$L^q$-extension from an $L^2$-$L^2$-result in Iosevich, Krause, Sawyer, Taylor and Uriarte-Tuero \cite{alex}, we state our second result. 

\begin{theorem}\label{second}
Let $\lambda$ be a tempered distribution whose Fourier transform is a locally integrable function satisfying the decay
$$|\hat\lambda(\cdot)|\leq C|\cdot|^{-\alpha}$$
and assume the condition \eqref{ball} on $\mu$. Then if $s>d-\alpha$, the operator 
$$Af:=\lambda*(f\mu)$$
is bounded from $L^p(\mu)$ to $L^2(\mu)$ for $\frac{1}{p}\in[\frac{1}{2},\frac{3s+2\alpha-2d}{2s})$, and from $L^p(\mu)$ to $L^{p'}(\mu)$ for $\frac{1}{p}\in[\frac{1}{2},\frac{2s+\alpha-d}{2s})$. The thresholds for $s$ and $p$ are sharp in general, exact values of endpoints aside.
\end{theorem}
\noindent
We remark that unlike Theorem \ref{main}, in this statement, we do not obtain a complete Riesz diagram but only the sharpness in the $q=p'$ and the $q=2$-direction. Like previously, the diagram will be symmetrical and one can always go vertically upward.
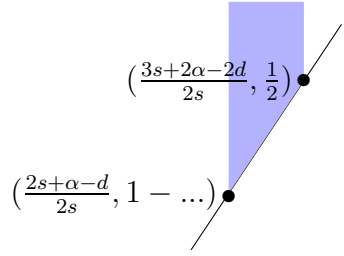
\begin{center}
\begin{tikzpicture}
\draw (0,0)--(2,3);
\path[fill=blue!30] (0.5,0.75)--(1.5,2.25)--(1.5,3.3)--(0.5,3.3)--(0.5,0.75);
\draw (0.5,0.75) node[left]{$(\frac{2s+\alpha-d}{2s},1-...)$};
\draw (0.5,0.7) node{$\bullet$};
\draw (1.5,2.25) node[left]{$(\frac{3s+2\alpha-2d}{2s},\frac{1}{2})$};
\draw (1.5,2.25) node{$\bullet$};
\end{tikzpicture}
\end{center}
\noindent
The sharpness for both $s$ and $p$ is obtained by the example $\lambda:=|\cdot|^{-(d-\alpha)}$. Its Fourier transform is $c_\alpha|\cdot|^{-\alpha}$, which is already known by Gel'fand \cite{gelfand}. In the original $L^2$-statement, it is noted that one can also consider $A$ of the form
$$(Af)(x)=\int K(x,y)(f\mu)(y)dy,$$ where $K$ is some reasonable function kernel such that $A$ is bounded from $L^2(\R^d)$ to $L^2_\alpha(\R^d)$, and for all $f$ in the Schwartz class,
\be\label{spp}
support(\widehat{Af})\subset C\cdot support(\hat f).
\ee
Moreover, only an approximate version of this condition is needed. The proof and conclusions will be of the same nature. 

\smallskip
In comparison, in Theorem \ref{main}, the fact that we integrate over a codimension-$1$ surface provides useful $L^1$-$L^\infty$-information. Moreover, the integral form provides an essential disjointness property, that is, Lemma \ref{dsjt} here, seen also for example in Eswarathasan, Iosevich, and Taylor \cite{alex3}, for Littlewood-Paley pieces with indices apart, which compensates the lack of the support condition \eqref{spp}. But Theorem \ref{main} and Theorem \ref{second} share the $L^p$-$L^2$-endpoint.

\smallskip
It is worth mentioning that if $\mu$ is the local Lebesgue measure on a fixed $n$-dimensional plane, then we are actually looking at the general Radon transform \eqref{radon} on $\R^k$ without the fractal setting. For example, if $d=3$ and $n=2$, then $s$ is exactly at its threshold, and the Riesz diagram is much nicer.

\smallskip
{\textbf{Questions.}} How to obtain the endpoint result on both $s$ and $p$? How to drop the condition $\mu[B(x,r)]\geq cr^s$, which backs up the $L^p$-$L^p$-endpoint, or does it provide a different Riesz diagram than the upper one alone? How does the exponent $\alpha$ in the assumed $L^2$-based Sobolev space estimate effect the Riesz diagram, or in higher dimension does large $\alpha$ still provide sharpness? We only use an $L^2$-averaging property, namely Lemma \ref{L2}, to connect to the fractal measure setting, are there other properties possibly provided by other reasonable conditions useful here?

\smallskip
{\textbf{Notation.}} In the following two sections, the expression $f\lesssim g$ will mean $f\leq C g$ for some constant $C$, $\gtrsim$ is similar, and $\approx$ refers to both the directions. The pairing $\langle f,g\rangle$ will mean the integral $\int f\bar g$, that is, the inner product between functions $f$ and $g$ with respect to the measure in the context.

\smallskip
{\textbf{Acknowledgment.}} The author would like to thank Prof. Alexander Iosevich for his generous help and advice.

\section{Proof of positive results}
Let $\rho_k:=\rho(\frac{\cdot}{2^k})$ be the dyadic partition of unity, where $\rho$ is a fixed smooth nonnegative radial function that vanishes outside $B(0,2)\backslash B(0,1)$, and $\rho_0:=1-\sum\limits_{k>0}\rho_k$. For any $f\in L^1(\mu)$, where $\mu$ is a compactly supported measure, let $(f\mu)_k$ denote the $k$-th Littlewood-Payley piece of $f\mu$, that is, the inverse Fourier transform of $\rho_k\widehat{f\mu}$, where $\widehat{f\mu}:=\int f(x)e^{-2\pi i\xi x}d\mu(x)$ shares the notation with the Fourier transform $\hat g$ for $g$ in the Schwartz class.

\subsection{$L^p$-$L^{p'}$-endpoint of Theorem \ref{main}}\label{pf1}
Consider the operator
\be\label{k}
A_kf:=A(f\mu)_k
\ee
By the definition of $A$, it suffices to bound the operator norm sum of $A_k$. We will apply the Riesz-Thorin interpolation theorem. To get the $L^2(\mu)$-$L^2(\mu)$ operator norm bound of $A_k$, by duality, for any $f$ and $g$ with $\|f\|_{L^2(\mu)},\|g\|_{L^2(\mu)}\leq 1$, we need to bound, by Plancherel, the quantity
$$\langle A_kf,g\rangle_\mu
=\langle A(f\mu)_k,g\rangle_\mu
=\langle \widehat{A(f\mu)_k},\widehat{g\mu}\rangle
=\sum\limits_j\langle\widehat{A(f\mu)_k},\widehat{(g\mu)_j}\rangle.$$
Now write
$$\sum\limits_j\langle\widehat{A(f\mu)_k},\widehat{(g\mu)_j}\rangle
=\sum\limits_{j:|k-j|\lesssim1}...+\sum\limits_{j:|k-j|\gtrsim1}...
=:I+II,$$
where the constant that separate the two cases of $|k-j|$ is to be determined. In some of the following lines, the case $j=0$ needs slight modification for the argument to go through, but for simplicity we omit it. By Cauchy-Schwarz,
$$|\langle\widehat{A(f\mu)_k},\widehat{(g\mu)_j}\rangle|
\leq \|\widehat{A(f\mu)_k}\|_{L^2[B(0,2^{j+1}\backslash B(0,2^j)]}\|\widehat{(g\mu)}\|_{L^2(B(0,2^{j+1}))},$$
where by the Sobolev bound on $A$, Plancherel, and Lemma \ref{L2},
\bast
\|\widehat{A(f\mu)_k}\|_{L^2[B(0,2^{j+1})\backslash B(0,2^j)]}
&\leq&2^{-j\alpha}\|\widehat{A(f\mu)_k}|\cdot|^\alpha\|_{L^2[B(0,2^{j+1})\backslash B(0,2^j)]}\\
&\lesssim&2^{-j\alpha}\|(f\mu)_k\|_{L^2(\R^d)}\\
&\lesssim&2^{-j\alpha}2^{k\frac{d-s}{2}}\|f\|_{L^2(\mu)},
\east
and $$\|\widehat{(g\mu)}\|_{L^2[B(0,2^{j+1})]}\lesssim2^{j\frac{d-s}{2}}\|g\|_{L^2(\mu)}.$$
So we estimate
$$|I|\lesssim2^{k(d-s-\alpha)}.$$
As for $II$, the error term, by Lemma \ref{dsjt} and Cauchy-Schwarz,
\bast
|II|&\lesssim&(\sum\limits_{j:j<k,|k-j|\gtrsim1}+\sum\limits_{j:j>k,|k-j|\gtrsim1})2^{-N\max(j,k)}\|f\|_{L^2(\mu)}\|g\|_{L^2(\mu)}\\
&\leq&k2^{-Nk}+\sum\limits_{j>k}2^{-Nj}\\
&\lesssim& 2^{-N'k},
\east
where $N$ and $N'$ can be arbitrarily large. We have shown $$\|A_k\|_{L^2(\mu)-L^2(\mu)}\lesssim 2^{k(d-s-\alpha)}$$
By Lemma \ref{1inf},
$$\|A_k\|_{L^1(\mu)-L^\infty(\mu)}\lesssim 2^k$$
Then the interpolation yields the desired result, with the assumption $s>d-\alpha$ making the geometric sum convergent.

%%%%%%%%%%%%%%%%%%%%
\subsection{$L^p$-$L^p$-endpoint of Theorem \ref{main}}
\noindent
By subsection \ref{pf1} ,
$$\|A_k\|_{L^2(\mu)-L^2(\mu)}\lesssim 2^{k(d-s-\alpha)},$$
with $A_k$ defined in \ref{k}. We will prove the restrictive weak $L^1$-$L^1$-estimate
\be\label{weak}
\|A_k(\chi_E)\|_{L^1(\mu)}\lesssim 2^{k(d-s)}\|\chi_E\|_{L^1(\mu)}.
\ee
Then, like in subsection \ref{pf1}, the interpolation, followed by the summation over $k$ of $\|A_k\|_{L^p(\mu)-L^p(\mu)}$  under the condition $s>d-\alpha$, gives the $L^p$-$L^p$-estimate up to the endpoint. If $\mu(E)=0$, the estimate is trivially true. Otherwise, let $\{B(w,\delta)\}_I$ be the covering for $E$ provided by Lemma \ref{coverlemma}, where the metric space is the support of $\mu$, and we may set each $\delta\lesssim2^{-k}$. Then,
\bea
\|A_k(\chi_E)\|_{L^1(\mu)}&\leq&
\int\int [|\widehat{\rho_k}|*(\chi_E\mu)](y)|\psi(x,y)|d\sigma_x(y)d\mu(x)\label{L1n}\\
&=&\int\int I_x(z)\chi_Ed\mu(z) d\mu(x)\nonumber\\
&\lesssim&\sum\limits_I\int\int I_x(z)\chi_{B(w,\delta)}(z)d\mu(z) d\mu(x)\nonumber,
\eea
where
$$J_x(z):=\int|\widehat{\rho_k}(z-y)||\psi(x,y)|d\sigma_x(y)$$
We dyadically decompose $J_x(z)$ into a sequence of regions indexed by nonnegative integer $l$, the $l$-th of which is $\approx2^{l-k}$ away from the submanifold defined locally by the zero set of $\Phi(x,\cdot)$. Let $J_{x,l}$ denote the restriction of $J_x$ to the $l$-th region. For $k>0$,
\be\label{N}
|\widehat{\rho_k}(z-y)|
=2^{kd}|\hat\rho[(z-y)2^k]|\lesssim 2^{kd}|(z-y)2^k|^{-N}
\ee
for arbitrarily large $N$, and for $k=0$, the right hand side is also a legitimate bound. As a smooth measure on a codimension-$1$ submanifold,
\be\label{mnfd}
[|\psi(x,\cdot)|\sigma_x][B(w,r)]\lesssim r^{d-1}
\ee
uniformly in $w$ and $r$. A further dyadic decomposition then, by \eqref{N} and \eqref{mnfd}, yields
\bea
J_{x,l}(z)
&=&\sum\limits_{j\geq l}\int_{y\in B(z,2^{j-k+1})\backslash B(z,2^{j-k})}|\widehat{\rho_k}(z-y)||\psi(x,y)|d\sigma_x(y)\label{Ixl}\\
&\lesssim&\sum\limits_{j\geq l}2^{kd}(2^j)^{-N}(2^{j-k+1})^{d-1}\nonumber\\
&\lesssim&2^{k-N'l}\nonumber,
\eea
where $N'$ can be arbitrarily large. The intersection between $B(w,\delta)$ and the support of $J_{x,l}$ is nonempty only if $w$ is in a $\lesssim2^{l-k}$-thick neighborhood of the submanifold defined by $\Phi$, while the set 
$$\{x:\Phi(x,y)=0\ for\ some\ y\in B(w,r)\},$$
defined only around the support of $\psi$, is contained in a $\lesssim r$-thick neighborhood of a codimension-$1$ submanifold with the constant uniform in $w$ and small $r$, which, covered by $\approx\frac{1}{r^{d-1}}$ many balls of radius $\approx r$, has $\mu$-measure $\approx r^{s-d+1}$. This, along with \eqref{Ixl}, leads to 
\be\label{intI}
\int\int J_{x,l}(z)\chi_{B(w,\delta)}(z)d\mu(z)d\mu(x)
\lesssim2^{k-N'l}\delta^s(2^{l-k})^{s-d+1}=2^{k(d-s)-N''l}\delta^s,
\ee
where $N'$ and $N''$ can be arbitrarily large. Summing over $l$ and $I$, \eqref{L1n}, \eqref{intI}, and the comparability $\sum\limits_I\delta^s\approx\mu(E)$ provided by Lemma \ref{coverlemma} now imply the desired restricted weak estimate \eqref{weak}.

%%%%%%%%%%%%%%%%%%%%

\subsection{$L^p$-$L^2$-endpoint of Theorem \ref{main}}
Recall $A_k$ defined in \ref{k}. Note that by subsection \ref{pf1}, 
\be\label{22}
\|A_k\|_{L^2(\mu)-L^2(\mu)}\lesssim 2^{k(d-s-\alpha)}
\ee
Now, to obtain an $L^1$-$L^2$ norm bound for $A_k$ by duality, consider for any $f,g$ with $\|f\|_{L^1(\mu)},\|g\|_{L^2(\mu)}\leq1$, by Plancherel,
$$
\langle A_kf,g\rangle_\mu=
\sum\limits_j\langle\widehat{A(f\mu)_k},\widehat{(g\mu)_j}\rangle=\sum\limits_{|j-k|\lesssim1}+\sum\limits_{|j-k|\gtrsim1}=:I+II.
$$
Again, by Plancherel,
$$
I=\sum\limits_{|j-k|\lesssim1}\langle\widehat{(f\mu)_k},\widehat{A^*{(g\mu)_j}}\rangle,
$$
where $A^*$ denotes the adjoint of $A$. Then, except that we only bound $\|\widehat{(f\mu)_k}\|_{L^2(\R^d)}$ by $\|f\|_{L^1(\mu)}\|1\|_{L^2[B(0,2^{j+1})]}$, the same argument as in subsection \ref{pf1} yields
\be\label{12}
\|A_k\|_{L^1(\mu)-L^2(\mu)}\lesssim 2^{k(d-\frac{s}{2}-\alpha)}
\ee
An interpolation between \eqref{22} and \eqref{12} gives the $L^p$-$L^2$-result up to the endpoint.

\subsection{$L^p$-$L^{p'}$-endpoint of Theorem \ref{second}}\label{lambda}
Define
\be\label{lambdak}
A_kf:=\lambda*(f\mu)_k=\lambda*\widehat{\rho_k}*f\mu.
\ee
The same argument as in subsection \ref{pf1} then yields
\be\label{222}
\|A_k\|_{L^2(\mu)-L^2(\mu)}\lesssim 2^{k(d-s-\alpha)}
\ee
For any $f\in L^1(\mu)$, by the Fourier decay of $\hat\lambda$,
$$|\widehat{A_kf}|=|\hat\lambda\rho_k\widehat{f\mu}|\lesssim2^{-k\alpha}\|f\|_{L^1(\mu)},$$
implying
\be\label{111}
\|A_k\|_{L^1(\mu)-L^\infty(\mu)}\lesssim 2^{k(d-\alpha)}.
\ee
Under the condition $s>d-\alpha$, an interpolation between \eqref{222} and \eqref{111}, before summing over $k$, gives the $L^p$-$L^{p'}$-result up to the endpoint.

\subsection{$L^p$-$L^2$-endpoint of Theorem \ref{second}}
Still, $A_k$ defined in \ref{lambdak}, by the same argument as in subsection \ref{pf1},
\be\label{2222}
\|A_k\|_{L^2(\mu)-L^2(\mu)}\lesssim 2^{k(d-s-\alpha)}.
\ee
To obtain the result by interpolation, we need an $L^1$-$L^2$ norm bound for $A_k$. By duality, consider, for any $f\in L^1(\mu)$ and $g\in L^2(\mu)$,
$$|\langle A_kf,g\rangle_\mu|
=|\langle\hat\lambda\rho_k\widehat{f\mu},\widehat{g\mu}\rangle|
\lesssim 2^{-k\alpha}\|f\|_{L^1(\mu)}\|1\|_{L^2[B(0,2^{k+1})]}2^{k\frac{d-s}{2}}\|g\|_{L^2(\mu)},$$
where we use Plancherel, Cauchy-Schwarz, and Lemma \ref{L2}. This implies
\be\label{1111}
\|A_k\|_{L^1(\mu)-L^2(\mu)}\lesssim 2^{k(d-\frac{s}{2}-\alpha)}
\ee
Like before, an interpolation between \eqref{2222} and $\eqref{1111}$ gives the $L^p$-$L^2$-result.
%%%%%%%%%%%%%%%%%%%%
\subsection{An $L^2$-averaging Lemma}
\begin{lemma}\label{L2}
Let $\mu$ be a compactly supported measure satisfying the ball size upper bound \eqref{ball}. Then
$$\|\widehat{f\mu}\|_{L^2[B(0,R)]}\lesssim R^{\frac{d-s}{2}}\|f\|_{L^2(\mu)}.$$
\end{lemma}
\noindent
The proof is nicely written in Wolff \cite{wolff} (Theorem 7.4) and uses the Schur's test. 

\subsection{Lemma of essential disjointness}\label{dsjt}
\begin{lemma}
Let $\mu$ be a compactly supported measure, and $A$ is the generalized Radon transform \eqref{radon}, where $\Phi$ satisfies the nonvanishing partial gradient condition in Theorem \ref{main}. Then, for $|j-k|\gtrsim1$,
$$|\langle A(f\mu)_k,(g\mu)_j\rangle|
\lesssim
2^{-N\max(j,k)}\|f\|_{L^1(\mu)}\|g\|_{L^1(\mu)},$$
where $N$ can be arbitrarily large, and lower bound for $|j-k|$ is some constant that does not depend on the choices of $f,g$.
\end{lemma}\noindent For arbitrarily small $\epsilon>0$, we consider the operator
\begin{equation}\label{thick}
(A_\epsilon h)(x):=\frac{1}{\epsilon}\int h(y)\phi(\frac{\Phi(x,y)}{\epsilon})\psi(x,y)d(y),
\end{equation}
where $\phi$ is the fixed function described above Theorem \ref{main}. Then it suffices to establish the same estimate for $A_\epsilon$ uniformly in $\epsilon$. Use Plancherel, and write $(f\mu)_k$ and $\phi$ as the inverse Fourier transforms of their respective Fourier transforms, to get
\begin{equation}\label{integral}
\langle A(f\mu)_k,(g\mu)_j\rangle
=\int\widehat{A(f\mu)_k}\widehat{g\mu}\rho_j
=\int\int\int\widehat{f\mu}(\eta)\widehat{g\mu}(\xi)\hat\phi(\epsilon s)I_{j,k}(\xi,\eta,s)d\xi d\eta ds,
\end{equation}
where
$$I_{j,k}(\xi,\eta,s)
:=\rho_k(\eta)\rho_j(\xi)\int\int \psi(x,y)e^{2\pi i[\xi x+\eta y+s\Phi (x,y)]}dxdy$$
The gradient over $(x,y)$ of the phase $\xi x+\eta y+s\Phi(x,y)$ is $$[\xi+s\nabla_x\Phi(x,y),\eta+s\nabla_y\Phi(x,y)].$$
We claim
\begin{equation}\label{claim}
|I_{j,k}(\xi,\eta,s)|
\lesssim\rho_k(\eta)\rho_j(\xi)[\max(|\xi|,|\eta|,|s|)]^{-N}
\end{equation}
for arbitrarily large $N$. We present the case $k-j\gtrsim1$ in detail, and the argument for the other case is symmetrical. Assume that both gradients have norm values on $[a,b]$ for some $a,b\in\R^+$, and set the constant in $k-j\gtrsim1$ such that $|\xi|\leq\frac{a}{4b}|\eta|$ on the support of $I_{j,k}$. There are three subcases. i) If $|s|\geq\frac{2}{a}|\eta|$, then $$|\eta+s\nabla_y(x,y)|\geq\frac{a}{2}|s|.$$ Integration by parts with respect to the $y$-variable $N$ times yields
\begin{equation}\label{ibp}|I_{j,k}(\xi,\eta,s)|
\lesssim \rho_k(\eta)\rho_j(\xi)|\eta+s\nabla_y(x,y)|^{-N}
\lesssim \rho_k(\eta)\rho_j(\xi)[\max(|\xi|,|\eta|,|s|)]^{-N}.
\end{equation}
ii) If $|s|\leq\frac{1}{2b}|\eta|$, then
$$|\eta+s\nabla_y\Phi(x,y)|\geq\frac{1}{2}|\eta|.$$
Integration by parts with respect to the $y$-variable $N$ times yields \eqref{ibp} as well. iii) If $\frac{1}{2b}|\eta|\leq s\leq \frac{2}{a}|\eta|$, then $$|\xi+s\nabla_x\Phi(x,y)|\geq-|\xi|+\frac{a}{2b}|\eta|\geq\frac{a}{4b}|\eta|.$$ Integration by parts with respect to the $x$-variable $N$ times yields 
$$|I_{j,k}(\xi,\eta,s)|
\lesssim \rho_k(\eta)\rho_j(\xi)|\xi+s\nabla_x(x,y)|^{-N}
\lesssim \rho_k(\eta)\rho_j(\xi)[\max(|\xi|,|\eta|,|s|)]^{-N}.$$
With claim \eqref{claim}, in \eqref{integral}, we can first integrate over $s$ to get
$$|\langle A(f\mu)_k,(g\mu)_j\rangle|
\lesssim\int\rho_k(\eta)\rho_j(\xi)|\widehat{f\mu}(\eta)||\widehat{g\mu}(\xi)|[\\max(|\xi|,|\eta|)]^{-N}$$
for arbitrarily large $N$. Now integrate over $\xi$ and $\eta$ to reach the desired conclusion.

%%%%%%%%%%%%%%%%%%%%
\subsection{An $L^1$-$L^\infty$-Lemma}\label{1inf}
\begin{lemma}
Let $\mu$ be a compactly supported measure, and $A$ is the generalized Radon transform \eqref{radon}, where $\nabla_y\Phi(x,y)$ is nonzero on the support of $\psi$. Then, for every $x$,
$$A(f\mu)_k\lesssim 2^k\|f\|_{L^1(\mu)}$$
\end{lemma}\noindent
Like in the proof of Lemma \ref{dsjt}, we replace $A$ by $A_\epsilon$ defined by \eqref{thick}, and seek the estimate uniformly in $\epsilon$. Still, by writing both $(f\mu)_k$ and $\phi$ as the inverse Fourier transforms of their respective Fourier transforms,
\bast
&&[A_\epsilon(f\mu)_k](x)\\
&=&\int_{|s|\lesssim2^k}\hat\phi(\epsilon s)\int(f\mu)_k(y)\psi(x,y)e^{2\pi is\Phi(x,y)}dyds
+\int_{|s|\gtrsim2^k}\hat\phi(\epsilon s)I_{x,k}(s)ds\\
&=:&I+II,
\east
where $$I_{x,k}(s):=\int\widehat{f\mu}(\eta)\rho_k(\eta)\int\psi(x,y)e^{2\pi i[s\Phi(x,y)+\eta y]}dyd\eta$$
First we use the estimate
\bast
\int|(f\mu)_k(y)|dy
&=&2^{kd}\int|[\hat\rho(2^k\cdot)*(f\mu)](y)|dy\\
&\leq&2^{kd}\int|f(z)|\int|\hat\rho[2^k(y-z)]|dyd\mu(z)\\
&=&\int|f(z)|\int|\hat\rho(y-2^kz)|dyd\mu(z)\\
&\lesssim&\int|f(z)| d\mu(z)\\
&=&\|f\|_{L^1(\mu)}
\east
to control
$$|I|\lesssim 2^k\|f\|_{L^1(\mu)}.$$
Note that the gradient over $y$ of the phase $s\Phi(x,y)+\eta y$ in $I_{x,k}(s)$ is
$$s\nabla_y\Phi(x,y)+\eta,$$
whose norm $\gtrsim |s|$ if $|s|\gtrsim2^k$ and $|\eta|\lesssim 2^k$. Then integration by parts $N$ times yields
$$|II|\lesssim
\int_{s\gtrsim2^k}|\hat\phi(\epsilon s)|\int|\widehat{f\mu}|\rho_k(\eta)[\max(|s|,|\eta|)]^{-N}d\eta ds
\lesssim\|f\|_{L^1(\mu)},$$
completing the proof.

\subsection{A ball cover Lemma}\label{cover}
\begin{lemma}\label{coverlemma}
Let $\mu$ be a finite Borel measure on a separable metric space satisfying
$$\mu[B(x,r)]\approx r^s$$
uniformly for all $x$ and $r$ smaller than some fixed constant. Then any Borel set $E$ with $\mu(E)>0$ enjoys a countable ball cover $\{B(x,\delta)\}_I$ such that
$$\mu(E)\approx \sum\limits_I\delta^s,$$
where the constants do not depend on the choice of $E$.
\end{lemma}\noindent
Borel measure on metric space is outer regular, so we can choose an open set $O$ with $\mu(O)\leq 2\mu(E)$. Start with a collection of balls covering $E$, where each ball is contained in $O$. Then, by the Vitali Covering Lemma, it has a countable subcollection $\{B(x,\delta)\}_I$ such that $\{B(x,5\delta)\}_I$ covers $E$, and balls in
$\{B(x,\delta)\}_I$ are pairwise disjoint. Then $$\mu(E)\leq\sum\limits_I\mu[B(x,5\delta)]\lesssim\sum\limits_I\delta^s\lesssim\sum\limits_I\mu[B(x,\delta)]\leq\mu(O)\lesssim\mu(E),$$ completing the proof. 

\section{Sharpness examples}
The consideration of the Cartesian product of two Cantor measures in terms of sharpness can be seen in Mattila \cite{mattila}, and the argument in this section will be similar.
\subsection{Diagram in Theorem \ref{main}}\label{dia}
Set $\R^d=\R^{k+1}\times\R^{d-k-1}$, and fix an interval $I\in \R$. Let $\sigma$ be the Cartesian product between the unit sphere measure on $\R^{k+1}$ and the Lebesgue measure on $I^{d-k-1}$. Let $Af:=\sigma*f$. Then $A$ is bounded from $L^2(\R^d)$ to $L^2_{\alpha}(\R^d)$ with $\alpha=\frac{k}{2}$, since, with $k$ non-vanishing principle curvatures, $\sigma$ has Fourier decay $\frac{k}{2}$. See, for example, Stein and Shakarchi \cite{stein} for more detail. Note that although we choose the spherical average, the argument for the generalized Radon transform is essentially the same. Let $A_\epsilon$ denote the convolutional operator whose kernel is
$$\frac{1}{\epsilon}\chi_{\{1-\frac{\epsilon}{2}<|\cdot|<1+\frac{\epsilon}{2}\}}$$
Set a measure $\mu$ such that, for any small $\epsilon$,
$$\mu[B(0,\epsilon)\times I^{d-k-1}]
\approx\epsilon^{s-(d-k-1)}$$
and
$$\mu(\{(x,y)\in \R^{k+1}\times I^{d-k-1}: 1<|x|<1+\epsilon\})
\approx\epsilon^{s-d+1},$$
via, for example, the Cartesian product between the Cantor measure and the Lebesgue measure. Denote
$$f_\epsilon:=\chi_{B(0,\epsilon)\times I^{d-k-1}}.$$
Note that, on the set
$$\{(x,y)\in \R^{k+1}\times I^{d-k-1}: 1<|x|<1+C\epsilon\},$$
the function $A_\epsilon f_\epsilon\gtrsim\frac{1}{\epsilon}\epsilon^{s-(d-k-1)}
\approx\epsilon^{s-(d-k-1)}$.
If $A$ is bounded from $L^p(\mu)$ to $L^q(\mu)$, then, necessarily,
$$\frac{1}{\epsilon}\epsilon^{s-(d-k-1)}(\epsilon^{s-d+1})^{\frac{1}{q}}
\lesssim\|A_\epsilon f_\epsilon\|_{L^q(\mu)}
\lesssim\|f_\epsilon\|_{L^p(\mu)}
\approx(\epsilon^{s-(d-k-1)})^{\frac{1}{p}}.$$
Sending $\epsilon$ to $0$ yields
$$\frac{1}{p}\leq\frac{s-d+k+\frac{1}{q}(s-d+1)}{s-d+k+1}.$$
With $k=2\alpha$ this gives the sloped line in the diagram. 

\smallskip
As for the vertical restricting line in the diagram, we still use the same $A$ and the corresponding $A_\epsilon$. Let $L_I$ denote the Lebesgue measure on $I$ and $C_I$ the cantor measure on $I$ with $C_I[B(x,\epsilon)]\approx\epsilon^{s-d+1}$ for $x$ on its support. Set $\mu$ to be the Cartesian product
$$C_I\times L_I^k\times L_I^{d-k-1}$$ and nearby just the Lebesgue measure on $\R^d$. Denote $$f_\epsilon:=\chi_{I\times B(0,\epsilon)\times I^{d-k-1}}.$$
Then on a region of measure $\gtrsim1$, $A_\epsilon f_\epsilon$ has value $\gtrsim \frac{1}{\epsilon}\epsilon^k$. If $A$ is bounded from $L^p(\mu)$ to $L^q(\mu)$, then, necessarily,
$$\frac{1}{\epsilon}\epsilon^k
\approx\|A_\epsilon f_\epsilon\|_{L^q(\mu)}
\lesssim\|f_\epsilon\|_{L^p(\mu)}
\approx (\epsilon^k)^\frac{1}{p}.$$
Sending $\epsilon$ to $0$ gives
$$\frac{1}{p}\leq\frac{s-d+k}{k}$$

\smallskip
But we need $s-d+1>0$, provided $s>d-\alpha=d-\frac{k}{2}$. This means $d-\frac{k}{2}\geq d-1$ i.e. $k\leq2$. Thus, in dimension $2$ and $3$, the spherical convolution operator gives the sharpness example, while in higher dimension, we only convolutes with the unit sphere in the first $2$ or $3$ coordinates.

\subsection{Threshold for $s$ in Theorem \ref{main}}
We show that $s\geq d-\alpha$ is in general necessary for any $L^p$-$L^q$-estimate to exist by providing an $L^\infty$-example. We use the same $A$ and $A_\epsilon$ as in subsection \ref{dia}. Let $I$ be a fixed interval of length $2$. Let $C_I$ denote the Cantor measure with $C_I[B(x,\epsilon)]\approx\epsilon^{s-d+1}$ for $x$ on its support and $L_I$ the Lebesgue measure on $I$. Set $$\mu:=C_I\times L_I^k\times L_I^{d-k-1}$$
Then, since one can fit a rectangle of size $\approx\epsilon^\frac{1}{2}\times\epsilon$ into an $\epsilon$-thick annulus, $A_\epsilon(1\mu)$ on a region with $\mu$-measure $\approx1$ has value
$$\frac{1}{\epsilon}(\epsilon^k)^\frac{1}{2}\epsilon^{s-d+1}.$$
Necessarily, this should be bounded, thus sending $\epsilon$ to $0$ gives
$$s\geq d-\frac{k}{2}=d-\alpha.$$
Note that, just as in subsection \ref{dia}, in dimension $4$ and higher we set $\alpha$ small, and if $s\in(d-l,d-l+1]$ for some integer $l>1$, we set the following. Let $C_I$ denote the Cantor measure on $I$ with $C_I[B(x,\epsilon)]\approx\epsilon^{s-(d-1)}$ for $x$ on its support, and let $\delta$ denote the Dirac measure at $0\in\R$. Set $\mu$ to be 
$$C_I\times\delta^{l-1}\times L_I^{d-l}$$
and the Lebesgue measure nearby. Then on a region with $\mu$-measure $\gtrsim1$,
$$A_\epsilon(1\mu)\gtrsim\frac{1}{\epsilon}(\epsilon^{s-d+l})^{\frac{1}{2}}.$$
Sending $\epsilon$ to $0$ gives $s\geq d-l+2$, a contradiction.

\smallskip

%%%%%%%%%%%%%%%%%%%%

\subsection{An example for Theorem \ref{second}}
Set $\lambda:=|\cdot|^{-(d-\alpha)}$. Then $\hat\lambda=c_\alpha|\cdot|^{-\alpha}$. Define, for any $f\in L^1(\mu)$,
$$A_\epsilon f:=\lambda*\phi_\epsilon*(f\mu).$$
Set $f_\epsilon:=\chi_{B(0,\epsilon)}$.
Then, on $B(0,\epsilon)$, $A_\epsilon f_\epsilon$ has value $\gtrsim \epsilon^s\epsilon^{-(d-\alpha)}=\epsilon^{s-d+\alpha}$. For $A$ to have $L^p(\mu)$-$L^q(\mu)$-bound, necessarily, 
$$\epsilon^{s-d+\alpha}(\epsilon^s)^{\frac{1}{q}}
\lesssim\|A_\epsilon f_\epsilon\|_{L^q(\mu)}
\lesssim\|f_\epsilon\|_{L^p(\mu)}
\lesssim (\epsilon^s)^{\frac{1}{p}}.$$
Sending $\epsilon$ to $0$ yields
$$\frac{1}{p}\leq\frac{s-d+\alpha+s\frac{1}{q}}{s},$$
which is the sloped restricting line in the diagram.
\smallskip

The sharpness for $s$ is examined in Iosevich, Krause, Sawyer, Taylor and Uriarte-Tuero \cite{alex} (Theorem 2.1, part ii).

\medskip

\end{document}